\documentclass[a4paper,12pt]{article}
\usepackage{amsmath,amsthm,amsfonts,amssymb,bm,mathrsfs}
\usepackage{enumerate}
\usepackage{graphicx}
\usepackage{natbib}
\usepackage{xcolor}
\usepackage[scale=0.75]{geometry}

\usepackage{setspace}
\newtheorem{defn}{Definition}
\newtheorem{thm}{Theorem}
\newtheorem{lem}{Lemma}
\newtheorem{prop}{Proposition}
\newtheorem{rmk}{Remark}
\newtheorem{fact}{Fact}
\newtheorem*{assum}{Assumption}

\newcommand{\bR}{\mathbb{R}}

\newcommand{\cB}{\mathcal{B}}

\newcommand{\cF}{\mathcal{F}}

\newcommand{\cI}{\mathcal{I}}

\newcommand{\cL}{\mathcal{L}}
\newcommand{\cM}{\mathcal{M}}

\newcommand{\cT}{\mathcal{T}}

\DeclareMathOperator{\rmd}{d\!}

\DeclareMathOperator*{\argmax}{arg\,max}

\DeclareMathOperator{\supp}{supp}

\usepackage{hyperref}
\hypersetup{colorlinks=true,
           citecolor=blue,
           linkcolor=blue,
           urlcolor=blue,
           bookmarksopen=true,
           pdfstartview={XYZ null null 1.00},
           pdfpagelayout={SinglePage}
}

\begin{document}
\title{On the diffuseness of incomplete information game\thanks{We thank Yeneng~Sun and Haomiao~Yu for encouragement and discussion. The first author also thanks Nicholas~C.~Yannelis for several discussions when he visited University of Iowa on January 2013.}}
\author{Wei~He\thanks{Department of Mathematics, National University of Singapore, 10 Lower Kent Ridge Road, Singapore 119076. E-mail: he.wei2126@gmail.com.}
\and
Xiang~Sun\thanks{Department of Mathematics, National University of Singapore, 10 Lower Kent Ridge Road, Singapore 119076. E-mail: xiangsun.econ@gmail.com.}
\thanks{Department of Mathematical Economics and Finance, Economics and Management School, Wuhan University, Wuhan, Hubei, 430072, China.}
}
\maketitle

\abstract{We introduce the ``relative diffuseness'' assumption to characterize the differences between payoff-relevant and strategy-relevant diffuseness of information. Based on this assumption, the existence of pure strategy equilibria in games with incomplete information and general action spaces can be obtained. Moreover, we introduce a new notion of ``undistinguishable purification'' which strengthens the standard purification concept, and its existence follows from the relative diffuseness assumption.

\smallskip
\textbf{JEL classification}: C62; C70; C72; D82

\smallskip
\textbf{Keywords}: Game with incomplete information; Pure strategy equilibrium; Relative diffuseness; Undistinguishable purification
}

\newpage
\tableofcontents
\newpage
\setlength{\parskip}{4pt}

\section{Introduction}\label{sec-intro}

Since \cite{Harsanyi1967-68}, games with incomplete information have been widely studied and found applications in many fields. Various kinds of hypotheses are proposed on the formulation of such games to guarantee the existence of pure strategy equilibria.\footnote{The increasing literature has widened significantly in recent years, as evidenced by \cite{Athey2001}, \cite{AC2009}, \cite{Reny2011}, \textit{etc}.} In particular, if players' information is diffuse,\footnote{The information is said to be diffuse if every private information space is atomless.} positive results have been obtained when all players' action spaces are finite and the information structure is disparate; see \cite{RR1982}. These results lead to a natural conjecture on the existence of pure strategy equilibria in games with incomplete information and general action spaces; for example, see Theorem~6.2 of \cite{FT1991}. Unfortunately, this existence result fails to hold with the general action spaces; see \cite{KRS1999}.


The main aim of this paper is to consider the existence of pure strategy equilibria for games with incomplete information and general action spaces. Towards this end, we shall distinguish different roles of the diffuseness of information. In games with incomplete information, the private information will influence games from two aspects: payoffs and strategies. In the conventional approach, the different diffuseness of information on these two aspects is usually considered from a unified point of view. However, when making decisions, the player's strategy-relevant diffuseness of information could be much richer than that conveyed in the payoff function.\footnote{A player may have richer diffuseness of information when making decisions: from a realistic point of view, she can access to more information via communication, learning, \textit{etc}; from a technical point of view, to guarantee the measurability of her payoff function, one may only need a sub-$\sigma$-algebra.}

Therefore, we suggest to describe the strategy-relevant and payoff-relevant diffuseness of information separately. The relation between these two kinds of diffuseness is characterized by the ``relative diffuseness'' assumption, which basically says that the strategy-relevant diffuseness is essentially richer than the payoff-relevant diffuseness on any nonnegligible information subset. Based on this assumption, we are able to prove the existence of pure strategy equilibria in games with incomplete information and general action spaces without invoking any existence result of mixed strategy equilibria.

To obtain the existence of pure strategy equilibria, the purification method is another powerful tool originated from \cite{DWW1951a}. In the case of finite actions, the purification method ensures that every mixed strategy equilibrium has a payoff equivalent pure strategy equilibrium when the private information is diffuse; see, for example, \cite{RR1982}, \cite{MW1985} and \cite{KRS2006}. Since the existence of mixed strategy equilibria has been established with great generality (see, for example, \cite{MW1985} and \cite{Fu2008}), a pure strategy equilibrium automatically exists. However, the above results strictly depend on the assumptions on action spaces (finite or countable). To obtain the purification results for general action spaces, the underlying information spaces have to be sufficiently rich; see \cite{LS2006} for Loeb spaces.

A pure strategy profile is said to be a purification of a mixed strategy profile if the expected payoffs/distributions of these two strategy profiles are the same for all players. However, when a player restricts her attention to a payoff-relevant information subset, this payoff equivalence may be not valid. In this paper we will introduce the notion of ``undistinguishable purification'' which retains the payoff/distribution equivalence universally on any payoff-relevant information subset for each player. Based on the relative diffuseness assumption, we will show that an undistinguishable purification exists for any mixed strategy profile. Consequently, the existence of pure strategy equilibria in games with incomplete information and general action spaces can be also obtained.

The rest of the paper is organized as follows. Section~\ref{sec-gii} presents the setup of games with incomplete information, and the key assumption---relative diffuseness is introduced in Section~\ref{sec-rd}. In Section~\ref{sec-exist}, we prove the existence of pure strategy equilibria for general action spaces. The notion of undistinguishable purification is introduced and discussed in Section~\ref{sec-puri}. Section~\ref{sec-conclusion} concludes the paper and the proof of Theorem~\ref{thm-sufficiency} is collected in Section~\ref{sec-appen}.

\section{Games with incomplete information}\label{sec-gii}

Games with incomplete information (henceforth games for short) can be described as follows. Each player $i$ observes an informational type $t_i$, whose values lie in some measurable space $(T_i,\cT_i)$. After observing the type, player $i$ selects an action $a_i$ from some compact metric space $A_i$ of feasible actions. We allow each player's payoff to depend on the actions chosen by all the players, and on her type as well. We define an information structure $\lambda$ for the game which is a joint probability on $T_1 \times T_2 \times \cdots \times T_n$.

To be precise, a game with incomplete information $\Gamma$ consists of five formal elements.
\begin{itemize}
\item The set of players: $I=\{1,2,\ldots,n\}$.
\item The set of actions available to each player: $\{A_i\}_{i\in I}$. Each $A_i$ is a compact metric space. Let $A=\times_{i=1}^n A_i$ denote the space of action profiles.
\item The (private) information space for each player: $\{T_i\}_{i\in I}$. Each $T_i$ is endowed with a $\sigma$-algebra $\cT_i$. Let $T=\times_{i=1}^nT_i$ and $\cT=\otimes_{i=1}^n\cT_i$.
\item The payoff functions: $\{u_i\}_{i\in I}$. Each $u_i$ is a mapping from $A\times T_i$ to $\bR$.
\item The information structure: $\lambda$, a probability measure on the measurable space $(T,\cT)$.
\end{itemize}

For payoff functions, we have the following standard assumption. Conditions (1) and (2) describe the measurability and continuity respectively, and Condition~(3) states an integrably bounded restriction.
\begin{assum}[P]
\label{assum-p}
For each $i\in I$, the payoff function $u_i$ satisfies the following requirements:
\begin{enumerate}[(1)]
\item $u_i$ is $\cB(A)\otimes\cT_i$-measurable on $A\times T_i$, where $\cB(A)$ is the Borel $\sigma$-algebra of $A$.
\item $u_i(\cdot,t_i)$ is continuous on $A$ for all $t_i\in T_i$.
\item There is a real-valued integrable function $h_i$ on $(T_i,\cT_i,\lambda_i)$, such that $|u_i|\le h_i$ on $A\times T_i$.
\end{enumerate}
\end{assum}

For each $i\in I$, associated with the information structure $\lambda$ is a marginal probability on each $T_i$ which we denote by $\lambda_i$. For these probabilities, we have the following assumption of independence.
\begin{assum}[I]\label{assum-i}
The private information of each player is independent of all other players' private information, \textit{i.e.}, $\lambda=\otimes_{i=1}^n\lambda_i$.
\end{assum}

This assumption can be weakened and correlations of information are allowed.\footnote{For detailed discussions, see Section~\ref{sec-conclusion}.} We adopt this basic setup for the sake of simplicity.

For each player $i\in I$, a mixed strategy (resp. pure strategy) is a measurable function from $T_i$ to $\cM(A_i)$ (resp. $A_i$). The set of all mixed strategies (resp. pure strategies) of player $i$ is denoted by $L_0^{\cT_i}\big(T_i,\cM(A_i)\big)$ (resp. $L_0^{\cT_i}(T_i,A_i)$). A pure strategy can be viewed as a mixed strategy by taking it as a Dirac measure for almost every $t_i$. As usual, we write $t_{-i}$ for an information profile of all players other than $i$, and $T_{-i}$ for the space of all such information profiles. We adopt similar notations for action profiles and strategy profiles.

Given a strategy profile $f=(f_1,f_2,\ldots,f_n)$ and a subset $E\in \cT_i$, the expected payoff of player $i$ on the event $E$ is
$$U_i^E(f)=\int_{E\times (\times_{j\neq i}T_j)} \int_{A} u_i(a_1,\ldots,a_n,t_i) \cdot \prod\nolimits_{j\in I}f_j(t_j,\rmd a_j)\rmd\lambda(t).$$
Taking $E=T_i$, $U_i^{T_i}(f)$ is exactly the expected payoff of player $i$, which is denoted by $U_i(f)$ for simplicity.

Given a strategy profile  $f$, let
$$V_i^f(a_i,t_i)=\int_{\times_{j\neq i}T_j} \int_{A_{-i}} u_i(a_1,\ldots,a_n,t_i) \cdot \prod\nolimits_{j\neq i}f_j(t_j,\rmd a_j) \cdot \prod\nolimits_{j\neq i} \rmd\lambda_j (t_j),$$
then $U_i^E(f)=\int_E \int_{A_i} V_i^f(a_i,t_i) f_i(t_i,\rmd a_i) \rmd \lambda_i(t_i)$  for any subset $E\in \cT_i$.

A mixed (resp. pure) strategy equilibrium is a mixed (resp. pure) strategy profile $f^*=(f^*_1,f_2^*,\ldots,f^*_n)$ such that $f^*_i$ maximizes $U_i(f_i,f^*_{-i})$ in $L_0^{\cT_i}\big(T_i,\cM(A_i)\big)$ (resp. $L_0^{\cT_i}(T_i,A_i)$) for each $i\in I$.

\section{Relative diffuseness of information}\label{sec-rd}

In this section, we will propose the ``relative diffuseness'' assumption as a characterization of the relation between two kinds of diffuseness of information.

For each $i\in I$, let $(T_i,\cT_i,\lambda_i)$ be an atomless probability space with a complete countably-additive probability measure $\lambda_i$. Let $\cF_i$ be a sub-$\sigma$-algebra of $\cT_i$. The $\sigma$-algebras $\cT_i$ and $\cF_i$ will represent the diffuseness of information from the aspect of strategies and from the aspect of payoffs respectively.\footnote{$\cF_i$ could be regarded as the sub-$\sigma$-algebra generated by the payoff $u_i$, in the sense that $\cF_i$ is the smallest $\sigma$-algebra which can guarantee the measurability of $u_i$.} The probability spaces $(T_i,\cT_i,\lambda_i)$ and $(T_i,\cF_i,\lambda_i)$ will be used to model the information space and the payoff-relevant information space respectively.

For any nonnegligible subset $D\in\cT_i$, the restricted probability space $(D,\cF_i^D,\lambda_i^D)$ is defined as follows: $\cF_i^D$ is the $\sigma$-algebra $\{ D\cap D' \colon D'\in\cF_i \}$ and $\lambda_i^D$ the probability measure re-scaled from the restriction of $\lambda_i$ to $\cF_i^D$. Furthermore, $(D,\cT_i^D,\lambda_i^D)$ can be defined similarly.

Now we are ready to present the condition of nowhere equivalence.\footnote{The condition of nowhere equivalence is introduced in \cite{HSS2013}. Proposition~1 there presents several equivalent statements of this condition.}
\begin{defn}
\label{defn-nowhere_equiv}
Following the notations above, $\cT_i$ is said to be \textbf{nowhere equivalent} to $\cF_i$ if for every $D\in\cT_i$ with $\lambda_i(D)>0$, there exists a $\cT_i$-measurable subset $D_0$ of $D$ such that $\lambda_i(D_0\triangle D_1)>0$ for any $D_1\in\cF_i^D$.
\end{defn}

The following statement is clear: if $(T_i,\cT_i,\lambda_i)$ is atomless and the sub-$\sigma$-algebra $\cF_i$ is purely atomic, then $\cT_i$ is nowhere equivalent to $\cF_i$.\footnote{Lemma~1 in \cite{HSS2013} shows that for any atomless measure space, there always exists an atomless sub-$\sigma$-algebra, to which the original $\sigma$-algebra is nowhere equivalent.}

We shall introduce the assumption of relative diffuseness as follows.
\begin{assum}[RD]
\label{assum-rd}
For each $i\in I$, $\cT_i$ is atomless and nowhere equivalent to $\cF_i$.
\end{assum}

This assumption implies that on any nonnegligible set $D\subseteq T_i$, $\cT_i^D$ is always essentially larger than $\cF_i^D$, which means that the strategy-relevant diffuseness of information is richer.


\section{Existence of pure strategy equilibria}\label{sec-exist}

In this section, we turn to the issue of the existence of pure strategy equilibria in games with incomplete information and general action spaces. First, we need to modify Assumption~(P) to make the statement consistent.

\begin{assum}[P']
\label{assum-p}
For each $i\in I$, suppose that $\cF_i$ is a countably-generated sub-$\sigma$-algebra of $\cT_i$.\footnote{It is natural to assume ``countably-generated'' condition. Since the range of $u_i$ is $\bR$ endowed with the Borel $\sigma$-algebra, the restricted $\sigma$-algebra on $T_i$ which is induced by $u_i$ is obviously countably generated.} The payoff function $u_i$ satisfies the following requirements:
\begin{enumerate}[(1)]
\item $u_i$ is $\cB(A)\otimes\cF_i$-measurable on $A\times T_i$, where $\cB(A)$ is the Borel $\sigma$-algebra of $A$.
\item $u_i(\cdot,t_i)$ is continuous on $A$ for all $t_i\in T_i$.
\item There is a real-valued integrable function $h_i$ on $(T_i,\cF_i,\lambda_i)$, such that $|u_i|\le h_i$ on $A\times T_i$.
\end{enumerate}
\end{assum}

The theorem on the existence of pure strategy equilibria is presented below.
We will provide a direct proof in the Section \ref{sec-appen} without invoking any results of mixed/behavioral strategies.\footnote{There are some discussions in the literature on the limitations and relations of mixed and behavioral strategies; see \cite{RR1982} and \cite{MW1985}. Our proof is free from these issues.}

\begin{thm}\label{thm-sufficiency}
In a game with incomplete information $\Gamma$, the action space $A_i$ is uncountable for each $i\in I$. If Assumptions~(P'), (I) and (RD) hold, then there exists a pure strategy equilibrium.
\end{thm}

There is a substantial literature on the existence of pure-strategy equilibria with general action spaces and more structures on the information space, \textit{e.g.}, Loeb spaces (see \cite{KS1999}) and saturated probability spaces (see \cite{KZ2012}). Among these conditions, the information spaces cannot contain any countably-generated part. Our result is more general in the sense that any atomless probability information space is allowed.

\begin{rmk}
Since an atomless probability space is saturated if and only if it is nowhere equivalent to any countably generated sub-$\sigma$-algebra (see \cite{HSS2013}), the existence of pure strategy equilibria in games with general action spaces and saturated information spaces is straightforward based on Theorem~\ref{thm-sufficiency}.
\end{rmk}

\section{Undistinguishable purification}\label{sec-puri}

In this section, several equivalence conditions for strategy profiles will be discussed. We shall introduce the notion of undistinguishable purification which refines the standard purification concept and prove its existence based on the relative diffuseness assumption. Then the existence of pure strategy equilibria follows directly.

\subsection{Universal equivalences}\label{subsec-equi}

In this subsection, we will propose the notions of ``universal payoff/distribution equivalence'' by requiring that each player cannot distinguish the expected payoffs/distributions on any payoff-relevant information subset.

\begin{defn}\label{defn-undistinguishable}
Let $f=(f_1,f_2,\ldots,f_n)$ and $g=(g_1,g_2,\ldots,g_n)$ be two mixed strategy profiles.
\begin{enumerate}
  \item For player $i$, the strategy profiles $f$ and $g$ are said to be distribution equivalent on the event $E\in\cF_i$ if $\int_{E} f_i(t_i,\cdot)\rmd \lambda_i(t_i) = \int_{E} g_i(t_i,\cdot)\rmd \lambda_i(t_i)$.  The strategy profiles $f$ and $g$ are said to be universally distribution equivalent if for each player $i \in I$, $\int_{E} f_i(t_i,\cdot)\rmd \lambda_i(t_i) = \int_{E} g_i(t_i,\cdot)\rmd \lambda_i(t_i)$\footnote{When $f_i$ is a pure strategy, $\int_Ef_i(t_i,B)\rmd\lambda_i(t_i)=\lambda_i\left(E\cap f_i^{-1}(B)\right)$ for any Borel subset $B\subseteq A_i$.} on every event $E\in \cF_i$.
  \item For player $i$, the strategy profiles $f$ and $g$ are said to be payoff equivalent on the event $E\in \cF_i$ if $U_i^E(f) = U_i^E(g)$.  The strategy profiles $f$ and $g$ are said to be universally payoff equivalent if for each player $i \in I$, $U_i^E(f) = U_i^E(g)$ on every event $E\in\cF_i$.
  \item Suppose that $f$ is a pure strategy profile. For player $i$, $f_i$ is said to be belief consistent with $g_i$ if $f_i(t_i)\in \supp g_i(t_i)$ for $\lambda_i$-almost all $t_i\in T_i$. Moreover, $f$ is said to be belief consistent with $g$ if they are belief consistent for each player $i\in I$.
\end{enumerate}
\end{defn}

The following result is clear, and its proof is given for completeness.

\begin{lem}\label{lem-puri}
Suppose that $(I,\cI,\eta)$ is an atomless probability space, $\cL$ is a countably generated sub-$\sigma$-algebra of $\cI$, and $X$ is a Polish space. For any two $\cI$-measurable measure valued mappings $h_1$ and $h_2$ from $I$ to $\cM(X)$, if  for any event $D\in\cL$,
$$\int_{D} h_1(i,\cdot)\rmd \eta(i) = \int_{D} h_2(i,\cdot)\rmd \eta(i),$$
then
$$\int_D\int_X \psi(i,x) h_1(i,\rmd x)\rmd \eta(i) = \int_D\int_X \psi(i,x) h_2(i,\rmd x)\rmd \eta(i)$$
for any $\cL\otimes \cB(X)$-measurable integrably bounded mapping $\psi$.\footnote{An $\cL\otimes\cB(X)$-measurable mapping $\psi$ is said to be integrably bounded if $|\psi(i,x)|\le \phi(i)$ on $X$ for some real-valued integrable function $\phi$ on $(I,\cL,\eta)$.}
\end{lem}

\begin{proof}
For any $D\in \cF$ and any Borel subset $B\subseteq X$, we have
$$
\int_{I}\int_{X} 1_D(i) 1_B(x) h_1(i,\rmd x) \rmd\eta(i) = \int_{I}\int_{X} 1_D(i) 1_B(x) h_2(i,\rmd x) \rmd\eta(i).
$$

Fix an $\cL\otimes \cB(X)$-measurable integrably bounded mapping $\psi$. Without loss of generality, we assume that it is nonnegative. Then $\psi$ is an increasing limit of a sequence of simple functions. Thus,
$$\int_D\int_X \psi(i,x) h_1(i,\rmd x)\rmd \eta(i) = \int_D\int_X \psi(i,x) h_2(i,\rmd x)\rmd \eta(i)$$
by the monotone convergence theorem.
\end{proof}

The following proposition shows that the universal payoff equivalence can be deduced by the universal distribution equivalence.

\begin{prop}\label{prop-dis-payoff}
Suppose that Assumptions~(I) and (P') hold. If mixed strategy profiles $f$ and $g$ are universally distribution equivalent, then they are universally payoff equivalent.
\end{prop}

\begin{proof}
Fix player $i\in I$. Given any $E\in \cF_i$, it suffices to show $U_i^E(f)=U_i^E(g)$.

First, we have
\begin{align*}
V_i^g(a_i,t_i)
& = \int_{A_{-i}} \int_{\times_{j\neq i}T_j} u_i(a_1,\ldots,a_n,t_i) \cdot \prod\nolimits_{j\neq i}g_j(t_j,\rmd a_j) \cdot \prod\nolimits_{j\neq i} \rmd\lambda_j (t_j)\\
& = \int_{A_{-i}} u_i(a_1,\ldots,a_n,t_i) \cdot \prod\nolimits_{j\neq i}\mu_j^{g_j}(\rmd a_j),
\end{align*}
where $\mu_j^{g_j}=\int_{T_j} g_j(t_j,\cdot)\lambda_j(t_j)$. Since $f$ and $g$ are universally distribution equivalent, $\mu_j^{f_j}=\mu_j^{g_j}$ for each $j\neq i$, which implies that $V^f_i(a_i,t_i) = V^g_i(a_i,t_i)$ for any $a_i\in A_i$ and $t_i\in T_i$.

Since $f$ and $g$ are universally distribution equivalent, for any $D\in \cF_i$ and any Borel subset $B\subseteq A_i$, we have
$$
\int_{T_i}\int_{A_i} 1_D(t_i) 1_B(a_i) f_i(t_i,\rmd a_i) \rmd\lambda_i(t_i) =
\int_{T_i}\int_{A_i} 1_D(t_i) 1_B(a_i) g_i(t_i,\rmd a_i) \rmd\lambda_i(t_i).
$$
Thus, by Lemma~\ref{lem-puri} we have
\begin{align*}
U_i^E(f)
& = \int_{E}\int_{A_i} V^f_i(a_i,t_i) f_i(t_i,\rmd a_i) \rmd\lambda_i(t_i)\\
& = \int_{E}\int_{A_i} V^g_i(a_i,t_i) g_i(t_i,\rmd a_i) \rmd\lambda_i(t_i) = U_i^E(g).
\end{align*}
The proof completes.
\end{proof}

\subsection{Undistinguishable purification}\label{subsec-undis}
In this subsection, we will introduce the notion of ``undistinguishable purification'' and prove its existence. Furthermore, based on the existence results of mixed strategy equilibria, the existence of pure strategy equilibria for general action spaces can be obtained by the purification method.

The following fact is shown in Lemma~4.4(iii) of \cite{HK1984}.
\begin{fact}\label{fact-dis}
Suppose that $(I,\cI,\eta)$ is an atomless probability space, $\cL$ is a countably generated sub-$\sigma$-algebra of $\cI$, and $X$ is a Polish space. If $\cI$ is nowhere equivalent to $\cL$,\footnote{In Lemma~4.4(iii) of \cite{HK1984}, $\cI$ is assumed to be conditional atomless over $\cL$, which is equivalent to our condition; for detailed discussions, see \cite{HSS2013}.} then for any $\cL$-measurable mapping $h_1$ from $I$ to $\cM(X)$, there exists an $\cI$-measurable mapping $h_2$ from $I$ to $X$, such that $\int_{D} h_1(i,\cdot)\rmd \eta(i) = \int_{D} h_2(i,\cdot)\rmd \eta(i)$ for any event $D\in\cL$.
\end{fact}

In the following we will present a general purification result, which is beyond the framework of games with incomplete information. The proof follows from Fact~\ref{fact-dis} and Lemma~\ref{lem-puri}.

\begin{prop}\label{prop-puri}
Under the conditions of Fact \ref{fact-dis}, for any $\cL$-measurable mapping $h_1$ from $I$ to $\cM(X)$, there exists an $\cI$-measurable mapping $h_2$ from $I$ to $X$, such that $$\int_{D} h_1(i,\cdot)\rmd \eta(i) = \int_{D} h_2(i,\cdot)\rmd \eta(i)$$ for any event $D\in\cL$, and
$$\int_D\int_X \psi(i,x) h_1(i,\rmd x)\rmd \eta(i) = \int_D \psi(i,h_2(i)) \rmd \eta(i)$$
for any $\cL\otimes \cB(X)$-measurable integrably bounded mapping $\psi$.
\end{prop}

\begin{rmk}\label{rem-general puri}
In the literature, the purification result is always stated as follow: given a measure-valued mapping from $I$ to $\cM(X)$ and at most countably many joint measurable integrably bounded functions from $I\times X$ to $\bR$, there exists a measurable mapping from $I$ to $X$ which yields the same expected payoffs/distributions.\footnote{See, for example, \cite{DWW1951a} and \cite{KRS2006}.}

Proposition~\ref{prop-puri} generalizes this result by showing that one can find a measurable mapping from $I$ to $X$ which has the same expected payoffs/distributions for any joint measurable integrably bounded function from $I\times X$ to $\bR$ and on any payoff-relevant subset.
\end{rmk}

Now we are ready to give the definitions of purification and undistinguishable purification.
\begin{defn}\label{defn-undis}
Suppose that $f$ is a pure strategy profile and $g$ is a mixed strategy profile.
\begin{enumerate}
  \item The pure strategy profile $f$ is said to be a purification of $g$ if they are distribution equivalent and payoff equivalent on $T_i$ for each $i\in I$, and belief consistent.
  \item The pure strategy profile $f$ is said to be an undistinguishable purification of $g$ if they are universally distribution equivalent, universally payoff equivalent, and belief consistent.
\end{enumerate}
\end{defn}

The first definition is standard, and the latter one strengthens the first one by requiring that the expected payoffs/distributions are the same on any payoff-relevant information subset.

In \cite{DWW1951a} and \cite{KRS2006}, the payoff equivalence, distribution equivalence and belief consistence are stated separately. In Proposition~\ref{prop-dis-payoff} we have proven that the universal distribution equivalence implies the universal payoff equivalence, and in the following we will show that the belief consistence is also implied by the universal distribution equivalence. Therefore, the notion of undistinguishable purification can be simplified such that it only depends on the universal distribution equivalence.

\begin{prop}\label{prop-dis-belief}
If a $\cT$-measurable pure strategy $f$ is universally distribution equivalent to an $\cF$-measurable mixed strategy $g$, then $f$ is belief consistent with $g$.
\end{prop}

\begin{proof}
Fix player $i$. Define a mapping $c$ from $T_i\times A_i$ to $\bR$ as $c(t_i,a_i)=1_{\supp g_i(t_i)}(a_i)$. Then $c$ is $\cF_i \otimes \cB(A_i)$-measurable.\footnote{This measurability is implied by the measurability of the correspondence $\supp g_i(t_i)$.} Since $f$ and $g$ are universally distribution equivalent, by Lemma \ref{lem-puri} we have
$$\int_{T_i} c(t_i,f_i(t_i)) \rmd\lambda_i(t_i) = \int_{T_i}\int_{A_i} c(t_i,a_i) g_i(t_i,\rmd a_i) \rmd\lambda_i(t_i) = 1.$$ Therefore, $c(t_i,f_i(t_i))= 1$ for $\lambda_i$-almost all $t_i\in T_i$, which implies that $f_i(t_i)\in \supp g_i(t_i)$ for $\lambda_i$-almost all $t_i\in T_i$.
\end{proof}

It follows from Fact~\ref{fact-dis} that any $\cF$-measurable mixed strategy profile\footnote{By $\cF$-measurable strategy profile $g$, we mean that $g_i$ is $\cF_i$-measurable for each $i$; similar interpretation applies for the $\cT$-measurable strategy profile $f$.} has a universally distribution equivalent $\cT$-measurable pure strategy profile. Thus, the existence of undistinguishable purifications is due to Propositions~\ref{prop-dis-payoff}, \ref{prop-dis-belief} and Fact~\ref{fact-dis}.
\begin{thm}\label{thm-undis}
Suppose that Assumptions (P'), (I) and (RD) hold.
Then every $\cF$-measurable mixed strategy profile has a $\cT$-measurable undistinguishable purification.
\end{thm}
\begin{rmk}\label{rmk-standard puri}
The extensions of the standard purification result as in \cite{DWW1951a} and \cite{KRS2006} to the general action spaces is clear in our setup. Moreover, since an atomless probability space is saturated if and only if it is nowhere equivalent to any countably generated sub-$\sigma$-algebra, the corresponding purification consequences with saturated probability spaces also follow from our result.
\end{rmk}

The existence of mixed strategy equilibria has been established with great generality; see, for example, \cite{MW1985} and \cite{Fu2008}. Thus, an undistinguishable purification can be obtained by Theorem~\ref{thm-undis} and easily shown to be a pure strategy equilibrium.

\section{Concluding remarks}\label{sec-conclusion}

In the current paper, we propose the relative diffuseness assumption to characterize the differences between payoff-relevant and strategy-relevant diffuseness of information. Based on this assumption, the existence of pure strategy equilibria in games with incomplete information and general action spaces is obtained. Moreover, we introduce the notion of undistinguishable purification and show its existence.


The model discussed in this paper is simple as the prior is assumed to be independent across all players' private information spaces. The relative diffuseness assumption can be used to deal with more general formulations. For example, we can introduce a new type of information which is payoff-relevant for all players. To be clear, suppose that $T_0=\{t_0^1,t_0^2,\ldots,t_0^m\}$ represents the space of common information and $\cT_0$ the power set of $T_0$; see \cite{MW1985}.
For each $i\in I$, player $i$'s payoff function $u_i$ will depend on the action profile, her own type, and the realized common information.
The information structure $\lambda$ will be a probability on $\big(T_0\times(\times_{i=1}^nT_i),\cT_0\otimes(\otimes_{i=1}^n\cT_i)\big)$.
We need the following conditional independence assumption: for each $j\in\{1,2,\ldots,m\}$, $\lambda^j=\otimes_{i=1}^n\lambda_i^j$, where $\lambda^j$ is the conditional probability measure on the product measurable space $(\times_{i=1}^nT_i,\otimes_{i=1}^n\cT_i)$ when the common information is $t_0^j$, and $\lambda_i^j$ is the marginal probability measure of $\lambda^j$ on $(T_i,\cT_i)$ for each $i\in I$. Then Assumption~(RD) will be restated under the probability measure $\lambda_i^j$ for each $i$ and $j$. With the conditional independence assumption, (P') and (RD), it can be easily checked that the main results on purification and existence of a pure strategy equilibrium still hold following the similar arguments.

More generalizations could be achieved along this line. In \cite{RR1982}, independent payoff-relevant and strategy-relevant private information is considered. \cite{FSYZ2007} introduced a new concept of strategy-relevant public information, and \cite{Fu2008} extended the space of public information to be countable and the set of common information to be a general probability space. In all of these papers, the purification results and existence of pure strategy equilibria are restricted in the case of finite actions. By appropriately adopting the relative diffuseness assumption and appealing to analogous arguments, our results can be naturally extended to their  settings with general action spaces.

\section{Appendix}\label{sec-appen}
Let $(T,\cT,\lambda)$ be a probability space and $X$ a topological space with the Borel $\sigma$-algebra $\cB(X)$. A correspondence from $T$ to $X$ is a mapping $F$ from $T$ to the family of nonempty subsets of $X$. A correspondence $F$ is said to be $\cT$-measurable if its graph belongs to the product $\sigma$-algebra $\cT\otimes\cB(X)$. A mapping $f\colon T\to X$ is called a selection of $F$ if $f(t)\in F(t)$ for $\lambda$-almost all $t\in T$; $f$ is said to be a $\cT$-measurable selection if it is $\cT$-measurable. The distribution of $F$ is defined as follows:
$$D_{F}^{\cT}=\{\lambda f^{-1} \mid f\text{ is a $\cT$-measurable selection of }F\}.
$$

A correspondence $F$ from $T$ to $X$ is said to be upper-hemicontinuous at $t_0\in T$ if for any open set $U\subseteq X$ that contains $F(t_0)$, there exists a neighborhood $V$ of $t_0$ such that for every $t\in V$, $F(t)\in U$. A correspondence is upper-hemicontinuous if it is upper-hemicontinuous at every point $t\in T$. Furthermore, let $\cM(X)$ denote the space of Borel probability measures on $X$ with the topology of weak convergence.

The following lemma is proved by \cite{HS2013b}, which is a preparation for the proof of the existence of pure-strategy equilibria in games with incomplete information.

\begin{lem}\label{lem-distribution}
Suppose that $(T,\cT,\lambda)$ is an atomless probability space and $\cT$ is nowhere equivalent to its countably generated sub-$\sigma$-algebra $\cF$. Then we have the following statements.
\begin{enumerate}
\item[A1.] For any compact-valued $\cF$-measurable correspondence $F$ from $T$ to $X$, $D_F^\cT$ is compact and convex.\footnote{To obtain the convexicity of $D_F^\cT$, it suffices to assume that $F$ is closed-valued.}
\item[A2.] Let $F$ be a compact-valued $\cF$-measurable correspondence from $(T,\cF,\lambda)$ to X, $Y$ a metric space, and $G$ a closed-valued correspondence from $T\times Y\to X$ such that
    \begin{enumerate}[(a)]
    \item for every $(t,y)\in T\times Y$, $G(t,y)\subseteq F(t)$;
    \item for every $y\in Y$, $G(\cdot,y)$ is $\cF$-measurable;
    \item for every $t\in T$, $G(t,\cdot)$ is upper-hemicontinuous;
    \end{enumerate}
    Then $H(y)=D_{G_y}^{\cT}$ is upper-hemicontinuous.
\end{enumerate}
\end{lem}

Now we are ready to prove Theorem~\ref{thm-sufficiency}.

\begin{proof}[Proof of Theorem~\ref{thm-sufficiency}]
To consider pure-strategy equilibria, we focus on the interim payoff function. For each $i\in I$, the interim payoff function $F_i$ from $T_i\times A_i\times\left(\times_{j=1}^{n}\cM(A_j)\right)$ to $\bR$ is defined as follows:
$$F_i(t_i,a_i,\tau_1,\ldots,\tau_n)=\int_{A_{-i}}u_i(a_i,a_{-i},t_i)\rmd \tau_{-i}(a_{-i}).
$$
It is clear that $F_i$ is continuous on $A_i\times\left(\times_{j=1}^{n}\cM(A_j)\right)$ and $\cF_i$-measurable on $T_i$. For each $i\in I$, the best response correspondence $G_i$ from $T_i\times\left(\times_{j=1}^{n}\cM(A_j)\right)$ to $A_i$ is given by
$$G_i(t_i,\tau_i,\ldots,\tau_n)=\argmax_{a_i\in A_i} F_i(t_i,a_i,\tau_i,\ldots,\tau_n).$$
For each $t_i$, Berge's maximal theorem implies that $G_i$ is nonempty, compact-valued, and upper-hemicontinuous on $\times_{j=1}^{n}\cM(A_j)$.

For any $(\tau_i,\ldots,\tau_n)$, $F_i$ is continuous on $A_i$ and $\cF_i$-measurable on $T_i$. Then Lemma~III.14 (page~70) and Lemma~III.39 (Application, page~86) of \cite{CV1977} asserts that $F_i$ is $\cF_i\otimes\cB(A_i)$-measurable and $G_i$ admits an $\cF_i$-measurable selection. Thus $D_{G_i(\cdot,\tau_1,\ldots,\tau_n)}^{\cT_i}$ is not empty. Moreover, it is convex, compact-valued, and upper-hemicontinuous on $\times_{j=1}^{n}\cM(A_j)$ by Lemma~\ref{lem-distribution}.\footnote{Note that $G_i(t_i,\tau_1,\ldots,\tau_n)\subseteq H_i(t_i)$, where $H_i$ is a correspondence from $T_i$ to $A_i$ such that $H_i(t_i)\equiv A_i$ for all $t_i$.}

Consider a correspondence $\psi$ from $\times_{j=1}^{n}\cM(A_j)$ to itself:
$$\psi(\tau_1,\ldots,\tau_n)=\times_{i=1}^{n}D^{\cT_i}_{G_i(\cdot,\tau_1,\ldots,\tau_n)}.$$
Then it is clear that $\psi$ is nonempty, convex, compact-valued, and upper-hemicontinuous on $\times_{i=1}^n\cM(A_j)$. By Fan-Glicksberg's fixed-point theorem, there exists a fixed point $(\tau^*_1,\ldots,\tau^*_n)$ of $\psi$. Thus for each $i$, there exists some $f_i^*\in L_0^{\cT_i}(T_i,A_i)$ such that $f_i^*$ is a selection of $G_i(\cdot,\tau^*_1,\ldots,\tau^*_n)$ and $\tau_i^*$ is induced by $f_i^*$.

So the payoff of player $i$ is
\begin{align*}
U_i(f^*)
& = \int_T u_i(f_i^*(t_i),f_{-i}^*(t_{-i}),t_i)\rmd \lambda\\
& = \int_{T_i}\int_{T_{-i}} u_i(f_i^*(t_i),f_{-i}^*(t_{-i}),t_i)\rmd \lambda_{-i}\rmd \lambda_{i}\\
& = \int_{T_{i}}\int_{A_{-i}} u_i(f_i^*(t_i),a_{-i},t_i)\rmd \tau^*_{-i}\rmd \lambda_{i}.
\end{align*}

The first equality holds due to the definition of $U_i$, second equality holds based on Assumption~(I), and the third equality relies on change of variables. By the choice of $\tau^*$, we have that $(f_1^*,f_2^*\ldots,f_n^*)$ is a pure-strategy equilibrium.
\end{proof}

\end{document}